\documentclass[12pt]{article}

\usepackage{amsmath,amssymb}

\usepackage{graphicx}
\usepackage{epsfig}

\newtheorem{dfn}{Definition}
\newtheorem{thm}{Theorem}
\newtheorem{lem}{Lemma}[thm]
\newtheorem{prp}{Proposition}
\newtheorem{subl}{Lemma}[lem]

\def\R{{\mathbb R}}

\def\Sph{{\mathbb S}}
\def\phi{\varphi}
\def\epsilon{\varepsilon}
\def\rk{\mathrm{rk}\,}

\title{On the hull numbers of torus links}
\author{Ivan Izmestiev}

\begin{document}
\maketitle

\setcounter{section}{-1}

\section{Introduction}

Introduced in \cite{CKKS} by Jason Cantarella, Greg Kuperberg, Robert B. Kusner, and John M. Sullivan, the hull number is a complexity measure for links similar to the bridge number. By definition, it is the maximum integer $n$ such that for any realization of the link there is a point $c\in\R^3$ such that any plane through $c$ intersects the curves representing the link at least $2n$ times. Clearly, the hull number is less than or equal to the bridge number.

An elegant statement proved in \cite{CKKS} says that the hull number of a non-trivial knot is at least~$2$. By an integral-geometric argument this yields a new proof of F\'ary-Milnor theorem. The authors posed the question of finding links and knots with large hull numbers. Unexpectedly enough, it is difficult already to find a knot with hull number at least $3$. It was conjectured that the $(3,4)$-torus knot has hull number $3$. However, we show this to be false (see Figure 1).

In the present paper we study hull numbers of torus links and prove two lower bounds. First, for non-trivial torus links with $p$ components (and more generally, for any link with $p$ pairwise non-trivially linked components) we show that the hull number is at least $\frac{3}{5}p$. This bound is sharp in the sense that for any $p$ there are $p$-component links with components non-trivially linked and with hull number $\lceil \frac{3}{5}p \rceil$. In particular, these are $(p,p)$-torus links (``multi-Hopf links''). Secondly, for the hull number of a $(p,q)$-torus link with $p \le q$ we prove the lower bound $\frac{1}{2}p$. In both cases, the scheme of the proof is to show that a plane intersecting the link in a small number of points cuts off an accordingly small portion of the link, and to apply Helly's theorem. The idea to use Helly's theorem for the study of hull numbers appeared already in \cite{CKKS}.

In Section 1 we give definitions and state the main results. In Section 2 we prove the bound $\frac{3}{5}p$. The bound $\frac{1}{2}p$ is proved in Section 3 modulo a Lemma whose proof is given in Section 4. In Section 5 we construct some examples of links with a small hull number. Finally, in Section 6 we suggest a notion of hull number of a polytope.

The author wishes to thank G\"unter Rote, John Sullivan, and G\"unter Ziegler for helpful discussions, and Mathias Schacht and Vojtech R\"odl for pointing out the work \cite{HS}.

\section{Definitions and statements}
A link is an ambient isotopy class of a union of disjoint simple closed curves in $\R^3$. By abuse of notation we will denote with the same letter $L$ a link and a representative of it. A knot is a one-component link.

Following \cite{CKKS}, we define the \emph{$n$-th hull} of a union of curves and the \emph{hull number} of a link.
\begin{dfn}
Let $L$ be a union of disjoint simple closed curves in $\R^3$. The $n$-th hull $h_n(L)$ of $L$ is the set of points $p\in\R^3$ such that
each plane through $p$ cuts $L$ at least $2n$ times:
$$
h_n(L) := \{ c \in \R^3 \mbox{ such that } |S \cap L| \ge 2n \mbox{ for any plane } S \ni c \}
$$
\end{dfn}
Sometimes we have to count an intersection point of $L$ with a plane $S$ twice. Namely, it holds for any point where the curve ``does not changes sides'' with respect to $S$. The notion of changing sides is well-defined if the point is an isolated point of intersection. If the latter is not the case, the total number of intersections is infinite and therefore it doesn't matter how do we count the point. This counting rule ensures that the $n$-th hull is closed.

We have obvious inclusions $h_n(L) \subset h_{n+1}(L)$.
\begin{dfn}
We say that a union of curves $L$ has $n$ hulls, if we have $h_n(L) \ne \emptyset$ and $h_{n+1} = \emptyset$. For a link $L$ the hull number $u(L)$ is defined as the minimum number of hulls its representative can have. In other words,
$$
u(L) := \max \{ n |\: h_n(L) \ne \emptyset \mbox{ for any union of curves } L \mbox{ representing the link } L \}
$$
\end{dfn}

The main result of \cite{CKKS} says that for any non-trivial knot the hull number is at least~$2$. Surely, there must exist knots with bigger hull numbers. The following speculation suggests a candidate. We have an obvious inequality
$$
u(L) \le b(L),
$$
where $b(L)$ denotes the bridge number of $L$. On the other hand, it seems impossible to get rid of the $b(L)$-th hull provided that $L$ is a simple knot. In particular, for a torus knot of type $(p,q)$ the conjectured value of the hull number is $\min(p,q)$. However, it turns to be false at least for the $(3,4)$-torus knot.

\begin{prp}
Hull number of the $(3,4)$-torus knots is $2$.
\end{prp}
\textbf{Proof }See Figure \ref{fig:(3,4)}. $\Box$

\begin{figure}[htbp]
  \begin{center}
    \begin{picture}(0,0)%
\includegraphics{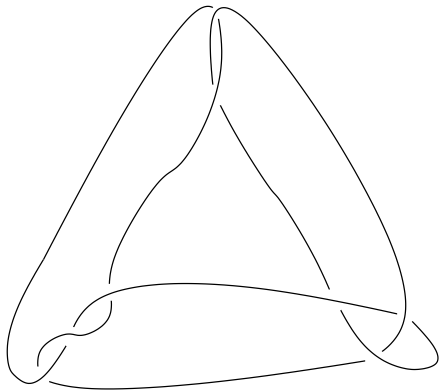}%
\end{picture}%
\setlength{\unitlength}{1036sp}%
\begingroup\makeatletter\ifx\SetFigFont\undefined
\def\x#1#2#3#4#5#6#7\relax{\def\x{#1#2#3#4#5#6}}%
\expandafter\x\fmtname xxxxxx\relax \def\y{splain}%
\ifx\x\y   
\gdef\SetFigFont#1#2#3{%
  \ifnum #1<17\tiny\else \ifnum #1<20\small\else
  \ifnum #1<24\normalsize\else \ifnum #1<29\large\else
  \ifnum #1<34\Large\else \ifnum #1<41\LARGE\else
     \huge\fi\fi\fi\fi\fi\fi
  \csname #3\endcsname}%
\else
\gdef\SetFigFont#1#2#3{\begingroup
  \count@#1\relax \ifnum 25<\count@\count@25\fi
  \def\x{\endgroup\@setsize\SetFigFont{#2pt}}%
  \expandafter\x
    \csname \romannumeral\the\count@ pt\expandafter\endcsname
    \csname @\romannumeral\the\count@ pt\endcsname
  \csname #3\endcsname}%
\fi
\fi\endgroup
\begin{picture}(10372,9227)(1079,-8596)
\end{picture}
    \caption{An incarnation of $(3,4)$-torus knot. Through any point there is a plane perpendicular to the diagram that intersects the knot
    at most $4$ times.}
    \label{fig:(3,4)}
  \end{center}
\end{figure}

In this paper we investigate the hull number of torus links and prove the following.

\begin{thm} \label{thm:TorLink}
For the $(p,q)$-torus link $L$ with $q$ divisible by $p$ (this link has $p$ trivial but non-trivially linked components) we have
$$
u(L) \ge \frac{3}{5} p
$$
\end{thm}

In general we have a weaker bound:
\begin{thm} \label{thm:Torus}
For the $(p,q)$-torus link $L$ with $p < q$ we have
$$
u(L) \ge \frac{1}{2} p
$$
\end{thm}

In Section \ref{sec:Upper} we also provide some exact values of $u(p,q)$. In particular, the lower bound from Theorem \ref{thm:TorLink} turns out to be sharp in the case $q=p$.

\section{Counting components}
In this section we prove a lower bound for the hull number of links with pairwise linked components. As a special case one obtains Theorem \ref{thm:TorLink}.

\begin{thm} \label{thm:Link}
Let $L = L_1 \sqcup \cdots \sqcup L_p$ be a $p$-component link in $\R^3$ with all of the components pairwise linked. Then there is a point $c \in \R^3$ such that each plane through $c$ meets at least $\frac{3}{5}p$ components of $L$. In particular, $u(L) \ge \frac{3}{5}p$.
\end{thm}

The proof is based on several Lemmas. First, some definitions.

\begin{dfn} \label{dfn:n-cut}
A plane $S$ is called an $n$-{\rm cut} (with respect to the link $L$) if it intersects at most $n$ components of $L$, all of them transversally.
\end{dfn}

\begin{dfn} \label{dfn:tiny}
A half-space is called {\rm tiny} (with respect to the link $L$) if it has an empty intersection with at least one of the components of $L$. The complement to a tiny half-space (thus, any half-space that contains one of the components in its interior) is called {\rm huge}.
\end{dfn}


\begin{lem} \label{lem:tiny}
Let $L$ be a link as in Theorem \ref{thm:Link} and $n < p$. Then any $n$-cut divides the space into a tiny and a huge half-spaces.
\end{lem}

\begin{lem} \label{lem:4tiny}
Let $L$ be a link as in Theorem \ref{thm:Link} and $n < \frac{3}{5}p$. Then $\R^3$ cannot be covered by four tiny half-spaces such that their boundaries are $n$-cuts with respect to $L$.
\end{lem}

Let us derive Theorem from the above Lemmas.

\noindent\textbf{Proof of Theorem \ref{thm:Link} }Take $n = \lceil \frac{3}{5}p \rceil - 1$, the largest integer less than $\frac{3}{5}p$. It follows from Lemma \ref{lem:4tiny} that any four huge half-spaces whose boundaries are $n$-cuts have a non-empty intersection. Helly's theorem \cite{DGK} implies that all huge half-spaces with boundaries $n$-cuts have a point in common.\footnote{One technical detail. Helly's theorem does not work in general for an infinite family of non-compact convex sets. One can proceed as follows. Take a closed ball $B \subset \R^3$ such that the set of curves $L$ lies in its interior. Lemma \ref{lem:4tiny} can be strengthened to say that $B$ cannot be covered as described. Then instead of huge half-spaces we can consider their (compact!) intersections with $B$.}
Let us show that any such point $c \in \R^3$ has the property in Theorem. Indeed, if it is not the case, then we have a plane $S \ni c$ that cuts at most $n$ components of $L$. We can slightly rotate $S$ to achieve transversality so that it becomes an $n$-cut passing through $c$. By Lemma \ref{lem:tiny} one of the half-spaces with boundary plane $S$ is huge. By shifting the plane $S$ a little so that it remains an $n$-cut we can achieve that this huge half-space does not contain the point $c$ any more. This contradicts the choice of $c$. $\Box$

\textbf{Proof of Lemma \ref{lem:tiny}} Denote the $n$-cut by $S$ and the associated half-spaces by $S_+$ and $S_-$. Note that $S_+$ and $S_-$ cannot each contain a whole component: otherwise, if $L_i \subset S_+$ and $L_j \subset S_-$, then $L_i$ and $L_j$ are unlinked. Thus at least one of the two half-spaces does not contain any whole component. Say, $S_-$ does not. Then it follows, firstly, that $S_+$ has a non-empty intersection with each of the link components, and secondly, that $S_-$ intersects at most $n$ of them, since $S_- \cap L_i \ne \emptyset$ implies $S \cap L_i \ne \emptyset$ and $S$ is an $n$-cut by assumption. Hence, $S_-$ is tiny and $S_+$ is huge. $\Box$

Lemma \ref{lem:4tiny}, in turn, follows from the two propositions below.

\begin{lem} \label{lem:AnyPair}
Let $L$ be a link as in Theorem \ref{thm:Link}. Suppose that we have a covering
$$
\R^3 = \bigcup_{\alpha=1}^m H_\alpha
$$
of $\R^3$ with half-spaces. Then for any two link components $L_i$ and $L_j$ there is an $\alpha$ such that $H_\alpha \cap L_i \ne \emptyset$ and $H_\alpha \cap L_j \ne \emptyset$.
\end{lem}

\begin{lem}[Hor\'ak, Sauer \cite{HS}] \label{lem:Graph}
Let the complete graph on $p$ vertices be covered by four complete subgraphs. Then the number of vertices in at least one of the covering subgraphs is greater or equal $\frac{3}{5}p$.
\end{lem}

\textbf{Proof of Lemma \ref{lem:4tiny}} Use reductio ad absurdum: assume $\R^3 = \cup_{\alpha=1}^4 H_\alpha$, where each $H_\alpha$ is tiny with the boundary $n$-cut. Then each of the sets
$$
I_\alpha := \{i |\: L_i \cap H_\alpha \ne \emptyset\}
$$
has at most $n$ elements. By Lemma \ref{lem:AnyPair}, the complete subgraphs spanned by subsets $I_\alpha$ cover the complete graph on the vertex set $\{1,\ldots,p\}$. Then by Lemma \ref{lem:Graph} there exists an $I_\alpha$ with $|I_\alpha| \ge \frac{3}{5}$ that contradicts the inequality $n < \frac{3}{5}$. $\Box$

\textbf{Proof of Lemma \ref{lem:AnyPair}} Assume the converse: for some $i$ and $j$, neither of the covering half-spaces intersects both components $L_i$ and $L_j$. Consider the region
$$
D := \bigcup_{H_\alpha \cap L_i \ne \emptyset} H_\alpha
$$
We have $D \supset L_i$ and $D \cap L_j = \emptyset$. But since $\R^3 \setminus D$ is convex, this means that $L_i$ and $L_j$ are not linked. Contradiction. $\Box$

Peter Hor\'ak and Norbert Sauer investigated in \cite{HS} the minimum constant $\lambda_k$ such that for any $p$ there is a covering of a complete graph on $p$ vertices with $k$ subgraphs of size $\lambda_k p + o(p)$. They found $\lambda_k$ for all $k \le 7$. In particular, $\lambda_4 = \frac{3}{5}$, which implies Lemma \ref{lem:Graph}. Note that this is better than a trivial bound $\lambda_4 \ge \frac{1}{2}$ obtained by counting the edges.

We would like to present here a short self-contained proof of Lemma \ref{lem:Graph}.

\textbf{Proof of Lemma \ref{lem:Graph}} Let $I$ be a set of cardinality $p$, and $I_1, I_2, I_3, I_4$ subsets of $I$ such that the following holds:
$$
\begin{array}{lc}
(*) & \mbox{for any } i,j \in I \mbox{ there is } \alpha \in \{1,2,3,4\} \mbox{ such that } i,j \in I_\alpha
\end{array}
$$
We have to show that for some $\alpha$ holds $|I_\alpha| \ge \frac{3}{5}p$.

Consider those pairs $(\alpha, \beta)$ that have the property $I_\alpha \cup I_\beta = I$. They can be encoded in a graph $\Gamma$ with
$$
\begin{array}{ll}
\mbox{the vertex set } & V = \{1,2,3,4\} \mbox{ and} \\
\mbox{the edge set } & E = \{(\alpha, \beta) |\: I_\alpha \cup I_\beta = I\}
\end{array}
$$
\textbf{Claim 1} If $(\alpha,\beta) \notin E$, then $(\gamma, \delta) \in E$, where different Greek letters stand for different vertices.
\par\noindent Proof: Otherwise there would exist an $i \notin I_\alpha \cup I_\beta$ and a $j \notin I_\gamma \cup I_\delta$ that violates the condition $(*)$.

\noindent \textbf{Claim 2} If $\Gamma$ contains a triangle, then $|I_\alpha| \ge \frac{2}{3}p$ for some $\alpha$.
\par\noindent Proof: Suppose that $(1,2),\, (2,3)$, and $(3,1)$ all are in $E$. Then it can easily be shown that each element of $I$ lies in at least two of the subsets $I_1,\, I_2,\, I_3$. The claim now follows by a simple counting argument.

\noindent \textbf{Claim 3} If $\Gamma$ does not contain a triangle, then $\Gamma$ is a star, that is $E = \{(\alpha, \beta), (\alpha, \gamma), (\alpha, \delta)\}$ for some $\alpha$.
\par\noindent Proof: It easily follows from Claim 1.

Thus it suffices to consider the case when the graph $\Gamma$ contains a star. Without loss of generality we can assume:
$$
I_1 \cup I_2 = I_1 \cup I_3 = I_1 \cup I_4 = I
$$

Note that if there is an element of $I$ that lies only in $I_1$, then condition~$(*)$ implies $I_1 = I$ and we are done.

Thus we may assume $I_1 \subset I_2 \cup I_3 \cup I_4$. Denote $|I_1| = k$. Then one of the sets $I_2$, $I_3$, $I_4$ contains at least one third of the elements of $I_1$. Let $|I_2 \cap I_1| \ge \frac{k}{3}$. Furthermore, since $I_1 \cup I_2 = I$, we have $|I_2 \setminus I_1| = p-k$. Therefore we have
$$
|I_2| = |I_2 \cap I_1| + |I_2 \setminus I_1| \ge \frac{k}{3} + (p - k) = p - \frac{2}{3}k
$$
If $k < \frac{3}{5}p$, then $p - \frac{2}{3}k > p - \frac{2}{3} \cdot \frac{3}{5}p = \frac{3}{5}p$. Thus either $|I_1| \ge \frac{3}{5}p$ or $|I_2| \ge \frac{3}{5}p$. $\Box$

\section{Counting meridians}
\setcounter{thm}{2}
\setcounter{lem}{0}
In this section $L$ denotes the torus link of type $(p,q)$ or a representing union of curves. We always assume $p \le q$. The group $\pi_1(\R^3 \setminus L)$ is usually denoted by $G$.

Proof of Theorem \ref{thm:Torus} is similar to that of Theorem \ref{thm:Link}. A major difference is that instead of counting the components we will now count the meridional generators of $G$. Lemmas \ref{lem:tiny} and \ref{lem:4tiny} don't much change in their appearence, but get completely different proofs. In contrast to Lemma \ref{lem:tiny}, the proof of its analogue requires some work and is postponed to the next section. In return, the analogue of Lemma \ref{lem:4tiny} has a simpler proof, due to a stronger assumption $n < \frac{1}{2}p$.

We will use results of Markus Rost and Heiner Zieschang \cite{RZ} on meridional generators of the group of a torus link.

\begin{dfn}
A meridional element is an element of $\pi_1(\R^3\setminus L)$ that can be represented by a meridian curve; a meridian
curve for $L$ is a curve that bounds a disk intersecting $L$ in exactly one point.
\end{dfn}

Rost and Zieschang proved in particular the following:
\begin{prp}[Rost, Zieschang \cite{RZ}] \label{prp:MerGen}
Let $G' \subset G$ be the subgroup generated by meridional elements $g_1,\ldots,g_k$ with $k < p$. Then $G'$ is a free group of rank at most $k$ (in particular, $G' \ne G$), and there is a free basis of $G'$ consisting of meridional elements.
\end{prp}

Now we have to modify the definitions from the previous section.

\begin{dfn} \label{dfn:n-points-cut}
A plane $S$ is called an $n$-{\rm point cut} (with respect to the link $L$) if it intersects $L$ transversally and in at most $n$ points.
\end{dfn}

\begin{dfn} \label{dfn:small}
A half-space $H$ is called {\rm small} (with respect to the link $L$) if the homomorphism $\pi_1(H\setminus L)\to \pi_1(\R^3\setminus L)$ is not surjective. Otherwise $H$ is called {\rm large}.
\end{dfn}

\begin{lem} \label{lem:small}
Let $S$ be an $n$-point cut with $n < p$. Then $S$ divides the space into a small and a large half-spaces. Besides, if $S_-$ is the small half-space, then the group $G_-$ is free of rank at most $n$.
\end{lem}
Clearly, the subgroup $G_- \subset G$ is generated by a finite set of meridional elements. Thus, in view of Proposition \ref{prp:MerGen}, group $G_-$ in Lemma \ref{lem:small} is \emph{freely generated} by at most $n$ meridional elements.
\begin{lem} \label{lem:4small}
Let $n <\frac{1}{2}p$. Then $\R^3$ cannot be covered by four small half-spaces such that their boundaries are $n$-point cuts.
\end{lem}

Theorem \ref{thm:Torus} follows from Lemmas \ref{lem:small} and \ref{lem:4small} by the same argument as in the proof of Theorem \ref{thm:Link}.

\textbf{Proof of Lemma \ref{lem:4small}}
Assume the converse. Let $S^{(1)},\, S^{(2)},\, S^{(3)},\, S^{(4)}$ be the boundary planes of the covering half-spaces. They can be assumed transversal to each other, as is not hard to show. Let $S^{(i)}_-,\, S^{(i)}_+$ denote the small, respectively the large half-space associated to $S^{(i)}$ for $i = 1,2,3,4$, and let $G^{(i)}_-,\, G^{(i)}_+ \subset G = \pi_1(\R^3 \setminus L)$ be the images of the fundamental group homomorphisms induced by inclusions $S^{(i)}_- \setminus L \to \R^3 \setminus L$ and $S^{(i)}_+ \setminus L \to \R^3 \setminus L$, respectively. Denote
$$
\begin{array}{rcl}
S^{(12)}_- & := & S^{(1)}_- \cup S^{(2)}_- \\
S^{(12)}_+ & := & S^{(1)}_+ \cap S^{(2)}_+ = \overline{\R^3 \setminus S^{(12)}_-}
\end{array}
$$
The regions $S^{(12)}_-$ and $S^{(12)}_+$ have common boundary that we denote by $S^{(12)}$. By a PL-transformation of $\R^3$ we can turn $S^{(12)}$ into a plane that is a $2n$-point cut (a non-transversal intersection with $L$ might happen on the line $S^{(1)} \cap S^{(2)}$, but we can assume this line disjoint with $L$ from the very beginning). Since $2n < p$, Lemma \ref{lem:small} applies and one of the groups $G^{(12)}_-, \, G^{(12)}_+$ has to be a proper subgroup of $G$. We claim that it is $G^{(12)}_-$. Indeed, this group is generated by $G^{(1)}_-$ and $G^{(2)}_-$ and consequently (remark after Lemma \ref{lem:small}) by at most $2n < p$ meridional elements (choose the base point on the line $S_1 \cap S_2$). Then by Proposition \ref{prp:MerGen} it cannot be equal to $G$. This implies also $G^{(12)}_+ = G$.

Now, the same is true with the indices $1$ and $2$ replaced by $3$ and $4$. In particular, $G^{(34)}_-$ is a proper subgroup of $G$. But from $\R^3 = \cup_{i=1}^4 S^{(i)}_-$ follows $S^{(12)}_+ \subset S^{(34)}_-$ and hence $G^{(12)}_+ \subset G^{(34)}_-$ which is a contradiction.\footnote{To be exact: It follows that $G^{(12)}_+$ is conjugate to a subgroup of $G^{(34)}_-$. Or, if one doesn't want to look for an appropriate place for the base point again and again, one can compactify $\R^3$ to $\Sph^3$ and choose the base point at infinity.} $\Box$

\section{Proof of Lemma \ref{lem:small}}
Let us compactify the space $\R^3$ by the point at infinity. Then a plane in $\R^3$ becomes a sphere in $\Sph^3$. The definition of an $n$-point cut carries over to an arbitrary sphere $S$, and it is sufficient to prove Lemma \ref{lem:small} for any spherical $n$-cut, replacing the word ``half-space'' by ``ball''.

Let $S_-,\, S_+$ be the balls associated to the $n$-point cut $S$. Denote by $G_-,\, G_+$ the corresponding subgroups in the link group $G$. Our task is to show that one of these subgroups coincides with $G$, while the other is free of rank at most $n$.

\setcounter{lem}{1}
\begin{subl} \label{lem:Nice}
At least one of the subgroups $G_-,\, G_+$ coincides with $G$.
\end{subl}
\textbf{Proof} There is an unknotted embedded torus $T \subset \R^3$ such that $L \subset T$. We may assume $T$ to be transversal to $S$. Let us show that each component $C$ of $T \cap S$ is a zero-homotopic curve on $T$. First, since $C \subset S$, it is an unknot in $\R^3$. Besides, $S$ provides $C$ with a trivial framing. Due to the transversality of $T$ and $S$, the embedding of $C$ in $T$ defines also a trivial framing on $C$. It follows that $C$ is either zero-homotopic or homotopic to a meridian or a parallel of $T$. If $C$ is homotopic to a meridian or parallel, then $|C \cap L| \ge p$. Besides, in this case it can easily be shown that there is another meridian, respectively another parallel, among the components of $T \cap S$. This contradicts the assumption $|L \cap S| < 2p$.

A collection of zero-homotopic curves on $T$ splits it into connected components so that one of these components is a torus with holes. If $S_+$ is the ball where this component lies, then we have $G_+ = G$. $\Box$

For convenience, let us extend the definition of rank of a free group by setting $\rk G := p$. Due to Proposition \ref{prp:MerGen} $\rk G'$ for a subgroup $G' \subset G$ generated by meridians is nothing else as the minimum number of meridional generators for $G'$.
\begin{subl} \label{lem:Rk}
Let $G^{(1)}, G^{(2)}$ be subgroups of $G$, both generated by meridional elements. Let $G'$ be the subgroup generated by $G^{(1)}$ and $G^{(2)}$. Then holds
$$
\rk G' \le \rk G^{(1)} + \rk G^{(2)}
$$
\end{subl}
\textbf{Proof} Let $\rk G^{(i)} = n_i$, $i=1,2$. If $n_1+n_2 \ge p$, then there is nothing to prove. So assume $n_1+n_2 < p$, in particular $n_1,n_2 < p$. Then by Proposition \ref{prp:MerGen} there are free bases of $G^{(1)}, G^{(2)}$ consisting of $n_1, n_2$ meridional elements respectively. It follows that the subgroup $G'$ is generated by $n_1+n_2 < p$ meridional elements. Then, again by Proposition \ref{prp:MerGen}, we have $\rk G' \le n_1+n_2$. $\Box$

Now assume Lemma \ref{lem:small} to be false and choose among all cuts $S$ with
\begin{equation} \label{eqn:Count_Ex}
G_+ = G \mbox{ and } \rk G_- > n := \frac{1}{2}|S \cap L|
\end{equation}
one with the smallest $n$. We will arrive to a contradiction by showing that this counterexample can be simplified to another one with a smaller $n$.

Consider again an unknotted torus $T$ containing $L$ and transversal to $S$. Each component $C_i$ of $T \cap S$ bounds a disk $W_i \subset T$. There is a partial order on the set $\{C_i\}$ defined by $C_i < C_j \Leftrightarrow W_i \subset W_j$. Let $C = \partial W$ be a minimum element with respect to this order relation. (If $T$ is disjoint with $S$, then $\rk G_- = 0 = n$.) Consider two cases.

\textbf{Case 1} $W \cap L \ne \emptyset$

In this case $W \cap L$ consists of a number of disjoint chords in $W$ (note that $L$ is transversal to $C$). The idea is to pull one of these chords out of $W$.

Choose a chord $l$ that cuts from $W$ a disk $W'$ that contains no other chords. Let $U(W')$ be a small open regular neighborhoood of $W'$ in $S_\epsilon$, where $S_\epsilon$ is that of the balls $S_+, S_-$ which contains $W$. Put
$$
S'_\epsilon := S_\epsilon \setminus U(W')
$$
This defines an $(n-1)$-cut $S'$ that divides $\Sph^3$ into balls $S'_\epsilon$ and $S'_{-\epsilon}$. (From the point of view of the sphere $S$ we pull the arc $l$ out of $S_\epsilon$.) Denote by $G'_\epsilon, G'_{-\epsilon}$ the corresponding subgroups of $G$ and look how much do they differ from $G_\epsilon, G_{-\epsilon}$.

It is easy to see that $G'_{-\epsilon} = G_{-\epsilon}$. On the other hand, the group $G_\epsilon$ as a subgroup of $G$ is generated by $G'_\epsilon$ and an element $g$ that is represented by a meridian running around $l$ inside $S_\epsilon$. It follows by Lemma \ref{lem:Rk} that $\rk G'_\epsilon \ge \rk G_\epsilon - 1$.

Note that, whether $\epsilon = +$ or $-$, $G'_+$ remains equal to $G$ since the intersection $T \cap S'_+$ is still vast enough. Thus we get an $(n-1)$-cut $S'$ with $G'_+ = G$ and $\rk G'_- > n-1$ that contradicts the choice of $S$.

\textbf{Case 2} $W \cap L = \emptyset$

If so, then $W$ subdivides one of the balls bounded by $S$ (let us again use the notation $S_\epsilon$ for this ball) into two balls $S_{\epsilon}^{(1)}$ and $S_{\epsilon}^{(2)}$. Consider two subcases.

First, let one of $S_{\epsilon}^{(1)}$, $S_{\epsilon}^{(2)}$, say $S_{\epsilon}^{(2)}$, have an empty intersection with $L$ (but not necessarily with $T$). Then we can replace $S_{\epsilon}$ by $S_{\epsilon}^{(1)} \setminus U(W)$, where $U(W)$ is a small open regular neighborhood of $W$. This transformation changes neither $G_+$ nor $G_-$. But as a result the number of components in $S \cap T$ gets smaller. Thus we can return to choosing a minimum component of $C$.

The second subcase is $S_\epsilon^{(1)} \cap L \ne \emptyset$, $S_{\epsilon}^{(2)} \cap L \ne \emptyset$. Let us show that in this situation replacing $S_\epsilon$ either by $S_{\epsilon}^{(1)} \setminus U(W)$ or by $S_{\epsilon}^{(2)} \setminus U(W)$ gives a counterexample with a smaller $n$. Denote by $G_\epsilon^{(1)}$, respectively by $G_{\epsilon}^{(2)}$, the subgroup of $G$ generated by meridians in $S_{\epsilon}^{(1)}$, respectively $S_{\epsilon}^{(2)}$. Then $G_\epsilon$ is generated by $G_{\epsilon}^{(1)}$ and $G_{\epsilon}^{(2)}$, and Lemma \ref{lem:Rk} implies
$$
\rk G_\epsilon^{(1)} + \rk G_\epsilon^{(2)} \ge \rk G_\epsilon > n
$$
On the other hand, if $|\partial S_\epsilon^{(1)} \cap L| = n_1$ and $|\partial S_\epsilon^{(2)} \cap L| = n_2$, then we have $n_1 + n_2 = n$. It follows that either $\rk G_\epsilon^{(1)} > n_1$ or $\rk G_\epsilon^{(2)} > n_2$. Moreover, the subgroup of $G$ generated by meridians in the complement to either of $S_{\epsilon}^{(1)}, S_{\epsilon}^{(2)}$ contains $G_{-\epsilon}$ and therefore has rank greater than $n$. Thus we have an $n_i$-point cut with ranks of subgroups on both sides greater than $n_i$. By Lemma \ref{lem:Nice} one of the subgroups coincides with $G$, and we arrive to a contradiction with minimality of $n$. Lemma \ref{lem:small} is proved.

\section{Upper bounds} \label{sec:Upper}
In order to estimate the hull number of a given link from above, one has to construct a representative of the link with a small number of non-empty hulls. In this section we provide constructions for some special cases and speculate a little on the general case.

First, let us show that the result of Theorem \ref{thm:Link} is the best that one can get by taking into account only the information about linking between the components. Namely, for every $p>0$ we construct a $p$-component link with components pairwise linked and the hull number $\lceil \frac{3}{5}p \rceil$. Consider a tetrahedron $ABCD$ in $\R^3$. Draw two simple disjoint curves each of that approximate the closed broken line $BCD$. Also draw for each of the edges $AB,\, AC,\, AD$ a simple closed curve running along the edge (once in each direction). Clearly, this can be done in such a manner that one gets $5$ pairwise linked disjoint closed curves. Besides, for any point in $\R^3$ there is a plane through it that cuts at most $3$ components of the link thus obtained: choose a face of the tetrahedron that is not too close to the given point and draw the plane parallel to this face. So, the link has only $3$ non-empty hulls. Figure \ref{fig:(5,5)} shows the $(5,5)$-torus link realized in the described way.

\begin{figure}[htbp]
  \begin{center}
    \begin{picture}(0,0)%
\includegraphics{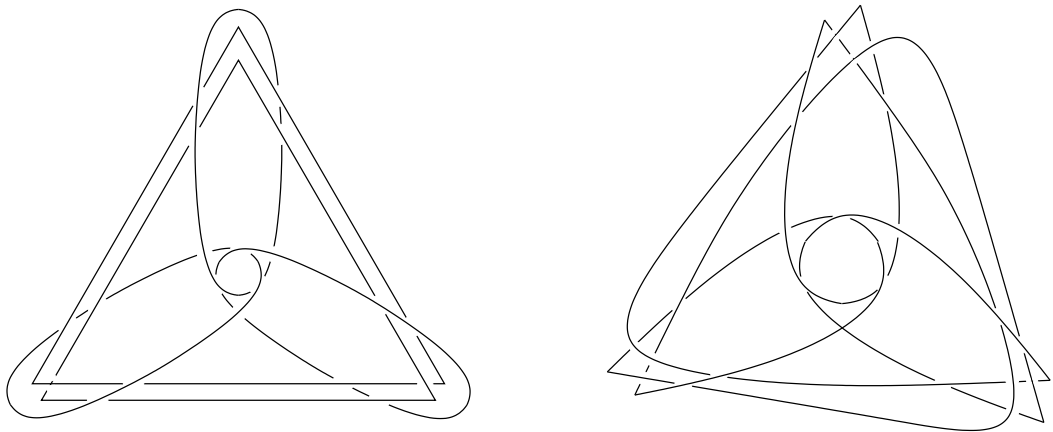}%
\end{picture}%
\setlength{\unitlength}{2072sp}%
\begingroup\makeatletter\ifx\SetFigFont\undefined%
\gdef\SetFigFont#1#2#3#4#5{%
  \reset@font\fontsize{#1}{#2pt}%
  \fontfamily{#3}\fontseries{#4}\fontshape{#5}%
  \selectfont}%
\fi\endgroup%
\begin{picture}(12524,5135)(524,-6106)
\end{picture}%

    \caption{Realizations of the $(5,5)$- and $(5,6)$-torus links with $3$ hulls. In each case, lift the central part of the diagram higher than the rest.}
    \label{fig:(5,5)}
  \end{center}
\end{figure}

Further, for any integer $p=5k+r,\, 0 \le r <5,$ one can similarly construct a link with $p$ components by putting $2k$ components along the circuit $BCD$ and $k$ components along each of the edges $AB,\, AC,\, AD$; if there is still something left ($r >0$), the remaining components should be distributed according to one of the patterns $(1,0,0,0),\, (2,0,0,0),\, (2,1,0,0),\, (2,1,1,0)$, where on the first place in each quadruple stands the number of components belonging to the triangular circuit. In the final picture one gets a link with $\lceil \frac{3}{5}p \rceil$ hulls. In particular, this construction can be realized for the $(p,p)$-torus link. One just has to cut the torus into $5$ strips that are thickenings of the components of the $(5,5)$-torus link, put the strips according to Figure \ref{fig:(5,5)} and draw on each of them an appropriate number of lines parallel to the boundary. 

Finally, the construction above can be generalized in the following way. Assign to each edge of a tetrahedron a non-negative integer, the multiplicity of an edge, in a way that the sum of multiplicities of the edges incident to any vertex would be even. Take $4$ balls centered at the vertices of a tetrahedron and join them by cylinders whose axes are the edges of the tetrahedron. Place into every cylinder a braid on the number of strands equal to the corresponding multiplicity, and into every ball the appropriate number of simple unknotted arcs (``bridges'') joining the ends of the strands. The braids and bridges can be isotoped so that the number of hulls of the link thus obtained equals half the maximum sum of multiplicities of edges incident to one vertex. As an example we give a diagram of the $(5,6)$-torus knot, see Figure \ref{fig:(5,5)}.

Thus we have the following proposition.
\begin{prp}
The hull number of the $(p,p)$-torus link is $\lceil \frac{3}{5}p \rceil$. The hull number of the $(5,6)$-torus knot is $3$.
\end{prp}

On the other hand, it seems that the estimate from Theorem \ref{thm:Link} is not sharp already for torus links. For example, we believe that $u(3,6) = 3$. In general we conjecture that for any $p$ there exists an integer $n$ such that the $(p,q)$-torus link with $q > n$ has hull number $p$.

\section{Final remarks}

The notion of the $n$-th hull allows a straightforward generalization from simple closed curves to graphs embedded in $\R^3$ (one should require that planes in the definition of the hull don't pass through the vertices of the graph; also, for convenience, each intersection point of a plane with the graph should be counted once, and not with the factor $\frac{1}{2}$). Then we define \emph{hull number of a polytope} as the minimum number of hulls of $1$-skeletons taken over all convex polytopes combinatorially equivalent to a given one. For example, hull number of the octahedron is $6$. Indeed, let $A_i^+,\, A_i^-,\, i = 1,2,3,$ be the vertices of a polytope, combinatorially equivalent to the octahedron, so that $A_i^+$ is opposite to $A_i^-$. Consider $3$ circuits $A_i^+ A_j^+ A_i^- A_j^-$, $i \ne j$. Their convex hulls always have a point in common, which thus belongs to the $6$-th hull of the polytope. We formulate this statement and its curious dualization in the last proposition of this paper.
\begin{prp}
Let $P \subset \R^3$ be a polytope combinatorially equivalent to the octahedron. Then there is a point in $\R^3$ such that each plane through this point intersects at least $6$ edges of $P$.

Let $Q \subset \R^3$ be any polytope combinatorially equivalent to the cube. Then there is a plane (usually at a long distance from the polytope, or even the plane at infinity) such that from each point of this plane at least $6$ vertices of the polytope are visible. 
\end{prp}

\end{document}